\documentclass[oneside,a4paper]{amsart} 
\usepackage{amssymb,verbatim,mathabx}

\usepackage[shortlabels]{enumitem}
    \setenumerate[0]{label={\rm (\roman*)}, leftmargin=*}

\usepackage[T1]{fontenc}
\usepackage[english]{babel}

\usepackage[textsize=footnotesize,textwidth=20ex,colorinlistoftodos]{todonotes}

\usepackage{hyperref}
\hypersetup{
    colorlinks=true,
    linkcolor=blue,
    filecolor=magenta, 
    citecolor=red,     
    urlcolor=cyan,
}

\usepackage[capitalize]{cleveref}
    \crefname{enumi}{}{}
    \Crefname{enumi}{Item}{Items}
    \crefname{equation}{}{}
    \Crefname{equation}{Equation}{Equations}

\newtheorem{proposition}{Proposition}[section]
\newtheorem{lemma}[proposition]{Lemma}
\newtheorem{corollary}[proposition]{Corollary}
\newtheorem{theorem}[proposition]{Theorem}

\theoremstyle{definition}

\newcommand\doverline[1]{\underline{#1}}

\begin{document}
\title{Whitney's Theorem for Line graphs of Multi-Graphs}
\author{Hans Cuypers} 
\maketitle

\begin{abstract}
Whitney's Theorem states that  every graph, different from
$K_3$ or $K_{1,3}$, is uniquely determined by its
line graph. 
A $1$-line graph of a multi-graph is the graph with as vertices the edges of the multi-graph, and two
edges adjacent if and only if there is a unique vertex on both edges.
The $\geq 1$-line graph of a multi-graph is the graph on the edges of the multi-graph, where 
two edges are adjacent if and only if there is at least one  vertex on both edges.
We extend Whitney's theorem to such line graphs of  multi-graphs, and show that most multi-graphs 
are uniquely determined by their line graph. 
Moreover, we present an algorithm to determine for a given graph $\Gamma$, if possible, a multi-graph
with $\Gamma$ as line graph.
\end{abstract}

\section{Introduction}

Let $\Delta=(V,E)$ be an ordinary graph, i.e. a graph without loops or multiple edges. So, we can consider each edge $e\in E$ as a set of two distinct vertices.
Then its line graph $L(\Delta)$ is the graph
with vertex set $E$ and two vertices $e,f\in E$ adjacent if and only if they intersect in a single vertex of $\Delta$. The graph $\Delta$ is called a \emph{root graph} for its line graph.

An \emph{isomorphism} $\phi$ between two graphs 
$\Delta=(V,E)$ and $\Delta'=(V',E')$ is a bijection $\phi:V\rightarrow V'$, such that
the map on $E$  induced by $\phi$, i.e. 
$\phi(\{v,w\})=\{\phi(v),\phi(w)\}$ for all $\{v,w\}\in E$, is a bijection from $E$ to $E'$.

Clearly, if  $\phi$ is an isomorphism between two graphs 
$\Delta=(V,E)$ and $\Delta'=(V',E')$, then
$\phi$ induces also an isomorphism between the two line graphs $L(\Delta)$ and $L(\Delta')$.

Whitney's Theorem is concerned with the reverse problem: 

\begin{itemize}
\item[]
Does an isomorphism between line graphs of graphs $\Delta$ and $\Delta'$ also imply the existence of an isomorphism between 
 $\Delta$ and $\Delta'$?
 
 \item[]In other words, is the root graph of a  line graph unique?
\end{itemize}
The answer to this question is stated in the following:

\begin{theorem}[Whitney's Theorem, \cite{whitney}]\label{whitney}
Suppose $\Delta$ and $\Delta'$ are connected graphs 
and $\phi$ an isomorphism between the line graph $L(\Delta)$ and $L(\Delta')$.

Then there is a unique isomorphism
$\psi$ between  $\Delta$ and $\Delta'$, inducing
$\phi$, unless
 $\Delta$ and $\Delta'$ are (up to permutation) $K_{1,3}$ and $K_3$.
\end{theorem}

An elementary and short proof of this result can be found in \cite{jung}.

In this note we consider the same problem for  multi-graphs.
Here a \emph{multi-graph} $\Delta=(V,E)$
consists of a set of vertices $V$ and a set of edges $E$ together with a relation $\sim$
between $V$ and $E$ such that
for any element $e\in E$ there are exactly two distinct elements $v,w\in V$ with $v,w\sim e$.
If $v\sim e$, then we will say $v$ is on $e$.
Of course, an ordinary graph is a multi-graph in which the relation $\sim$ is the relation $\in$.

There are at least two different obvious choices for the definition of the line graph of a multi-graph.
The vertices of the line graph of $\Delta=(V,E)$
are the elements from $E$. But two elements $e,f\in E$ can be defined to be adjacent
if and only if there is a unique vertex  with $v\sim e$ and $v\sim f$, or, 
if and only if there is at least one vertex $v$ with  $v\sim e$ and $v\sim f$.
We consider both situations. 

The  \emph{$1$-line graph} $L_1(\Delta)$ of a multi-graph $\Delta$ is the graph with vertex set $E$, such that two vertices $e,f\in E$ are adjacent 
if and only if there is a \emph{unique} vertex $v$ of $\Delta$ with $v\sim e$ and $v\sim f$.
By  $L_{\geq 1}(\Delta)$ we denote the \emph{$\geq 1$-line graph} of $\Delta$,
the graph with vertex set $E$ and two edges $e,f$ adjacent if there is a \emph{at least one} vertex $v$ of $\Delta$ with $v\sim e$ and $v\sim f$.

If it is clear from the context whether we are considering the $1$-line graph or $\geq 1$-line graph,
we will just refer to these graphs as line graphs.
In both cases  the graph $\Delta$ is called a \emph{root graph} for its line graph.

An \emph{isomorphism between two multi-graphs}
$\Delta=(V,E)$ and $\Delta'=(V',E')$
is a pair $(\phi_V,\phi_E)$
of bijections

$$\phi_V:V\rightarrow V', \phi_E:E\rightarrow E'$$
such that for all $v\in V$ and $e\in E$ we have
$$v\sim e\Leftrightarrow \phi_V(v)\sim\phi_E(e).$$ 

An isomorphism between two  multi-graphs
$\Delta=(V,E)$ and $\Delta'=(V',E')$ induces
an isomorphism, the isomorphism $\phi_E$, between the graphs $L_1(\Delta)$
and $L_1(\Delta')$, as well as between $L_{\geq 1}(\Delta)$
and $L_{\geq 1}(\Delta')$.

Again we can ask the following questions:
\begin{itemize}
\item[]
Does an isomorphism between line graphs of multi-graphs  $\Delta$ and $\Delta'$ also imply the existence of an isomorphism between 
 $\Delta$ and $\Delta'$?
 
 \item[]In other words, is the root graph of a  line graph of a multi-graph unique?
\end{itemize}

The answers to these questions are given in the following theorems:

\begin{theorem}\label{mainthm}
Suppose $\Delta$ and $\Delta'$ are connected multi-graphs 
and $\phi$ an isomorphism between the line graph $L_1(\Delta)$ and $L_1(\Delta')$.

Then there is a unique isomorphism
 between  $\Delta$ and $\Delta'$, inducing
$\phi$, unless
 $\Delta$ or $\Delta'$ contains four vertices.
 
 Moreover,  up to isomorphism, there is a unique connected multi-graph  $\Delta''$, not on $4$ vertices,
 such that $L_1(\Delta'')$ is isomorphic to $L_1(\Delta)$ and $L_1(\Delta')$.
\end{theorem}

\begin{theorem}\label{mainthm2}
Suppose $\Delta$ and $\Delta'$ are  connected multi-graphs 
and $\phi$ an isomorphism between the line graph $L_{\geq 1}(\Delta)$ and $L_{\geq 1}(\Delta')$.

Then there is a unique isomorphism
 between  $\Delta$ and $\Delta'$, inducing
$\phi$, unless
 $\Delta$ or $\Delta'$ contains a subgraph $\Delta_0$
 
 
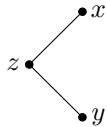
\begin{center}
 \begin{tikzpicture}[scale=0.7]
 \filldraw[color=black]
(0,0) circle [radius=2pt]
(1,1) circle [radius=2pt]
(1,-1) circle [radius=2pt];

\draw (1,1)--(0,0)--(1,-1);

\draw (1.3,1) node {$x$};
\draw (1.3,-1) node {$y$};
\draw (-0.3,0) node {$z$};

 \end{tikzpicture}
 \end{center}
 such that the vertices $x$ and $y$ have no other neighbors outside $\Delta_0$.
 
 Moreover, up to isomorphism, there is a unique connected multi-graph $\Delta''$ not containing such subgraph
 $\Delta_0$, with $L_{\geq 1}(\Delta'')$ isomorphic to $L_{\geq 1}(\Delta)$ and $L_{\geq 1}(\Delta')$.
\end{theorem}

Notice that $\Delta_0$ is not necessarily  an induced graph, we allow multiple edges on $x,z$ or $y,z$ or $x,y$.

After finishing the research for this note, we found that results similar to \cref{mainthm2} have also been obtained by
Zverovich \cite{analogue}.

 We present a uniform proof for the above theorems in the next section. This proof can be transformed into an algorithm that on input of an ordinary graph $\Gamma=(V,E)$ decides whether $\Gamma$ is a line graph of a multi-graph and also determines in the case that $\Gamma$ is indeed a line graph a root graph $\Delta$ for $\Gamma$.
 
The algorithm has complexity $O(|V|+|E|)$. 
It is discussed in \cref{sect:alg}.

A class of finite $1$-line graphs of multi-graphs that received a lot of attention is the class of so-called generalized line graphs, as they are graphs with  eigenvalues $\geq -2$. See for example
\cite{css} or \cite{book}, and various references in \cite{book}.

A \emph{generalized line graph} $\Gamma$ is a graph which can be constructed from a line graph of an ordinary graph $\Delta=(V,E)$ and  for each vertex
$v$ a cocktail party graph $\Delta_v$ (a possibly empty ordinary graph in which each vertex is
adjacent to precisely all but one other vertices).
The graph $\Gamma$ 
 is the union of the line graph
$L(\Delta)$ and all the graphs $\Delta_v$,
to which the following edges are  added:
a vertex $e$ of $L(\Delta)$ is adjacent to 
all vertices of $\Delta_v$, where $v$ is a vertex on $e$.

An equivalent definition is the following.
A graph $\Gamma$ is a generalized line graph if and only if it is the $1$-line graph of a 
multi-graph $\Delta$ such that:
\begin{enumerate}
\item Two vertices of $\Delta$ are on at most two common edges.
\item If $v,w$ are two vertices of $\Delta$ which are on two common edges, then
one of the vertices is on no other edges.
\end{enumerate}
See \cite[Proposition 6.2]{rootgraph}.

Using this definition of  a generalized line graph we can apply \cref{mainthm}
and easily deduce the analogue of Whitney's Theorem for generalized line graphs as in \cite{cameron} and \cite[Theorem 2.3.3, 2.3.4]{book}.

Moreover, our algorithm of \cref{sect:alg} can be modified for generalized line graphs. See also \cite{simic} for a different algorithm to find the root graph of a  generalized line graph.

\section{Line graphs}
\label{sect:proof}

In this section we prove \cref{mainthm} and  \cref{mainthm2}. 

We start with some definitions.

Let $\Gamma=(V,E)$ be an ordinary graph.
On the vertex set of $\Gamma$ we define
two equivalence relations $\equiv$ and $\bowtie$ by the following: two vertices $v,w$ of $\Gamma$ are in relation $\equiv$ if and only if they have the same set of neighbors (and hence are themselves not adjacent), while two vertices $v,w$ of $\Gamma$ are in relation $\bowtie$ if and only if they are adjacent and each other vertex adjacent to one of them is also adjacent to the other. 
For $v$ a vertex of $\Gamma$,
denote by  $[v]$ the $\equiv$-class of $v$ and by $\langle v\rangle$ the $\bowtie$-class of $v$.
Two $\equiv$-equivalent vertices are also called \emph{false twins}, while $\bowtie$-equivalent vertices are called \emph{true twins}.

By $\overline{\Gamma}$ we denote the graph with vertex set the $\equiv$-classes of  vertices of $\Gamma$, and two such classes $[v]$ and $[w]$  adjacent if and only if $v$ and $w$ are adjacent in $\Gamma$.
The graph on the $\bowtie$-classes, where  two such classes $\langle v\rangle$ and $\langle w\rangle$ are adjacent if and only if $v$ and $w$ are adjacent in $\Gamma$, is denoted by $\widehat{\Gamma}$.

If an $\equiv$- or $\bowtie$-class has size one, we often identify it with the element it contains.
In this way we justify a statement like $\overline{\Gamma}=\Gamma$, or  $\widehat{\Gamma}=\Gamma$, in case all equivalence classes have size one.

\begin{lemma}\label{except4}
Suppose $\Gamma$ is a line graph of a  connected ordinary  graph $\Delta$. Then the following hold:
\begin{enumerate}
\item 
$\overline\Gamma=\Gamma$ or $\Delta$ contains $4$ vertices.
\item
$\widehat\Gamma=\Gamma$ or $\Delta$ contains a (not necessarily induced)
subgraph $\Delta_0$
 
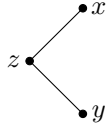
\begin{center}
 \begin{tikzpicture}[scale=0.7]
 \filldraw[color=black]
(0,0) circle [radius=2pt]
(1,1) circle [radius=2pt]
(1,-1) circle [radius=2pt];

\draw (1,1)--(0,0)--(1,-1);

\draw (1.3,1) node {$x$};
\draw (1.3,-1) node {$y$};
\draw (-0.3,0) node {$z$};

 \end{tikzpicture}
 \end{center}
 such that the vertices $x$ and $y$ have no other neighbors outside $\Delta_0$.
\end{enumerate}

\end{lemma}

\begin{proof}
Let $\Gamma=L(\Delta)$ for some graph $\Delta$.

Suppose $\overline\Gamma\neq\Gamma$.
Then there are two edges in $\Delta$, say $e\neq f$ that do not meet each other, such that each edge meeting one of $e$ or $f$ meets the other.
But then, by connectedness of $\Delta$, every edge of  $\Delta$ meets both $e$ and $f$ and we find $\Delta$ to be a connected graph on $4$ vertices.

Now  suppose  $\widehat\Gamma\neq\Gamma$.
Then there are two edges $e=\{x,z\}$ and $f=\{y,z\}$ of $\Delta$ on a common vertex $z$ such that each edge which meets $e$ also meets $f$.
But that implies that an edge meeting $e$ or $f$ meets them in $z$,
or it is the unique edge $\{x,y\}$.
So, the subgraph $\Delta_0$ of $\Delta$ on the vertices $x,y$ and $z$ and the edges $e$ and $f$ is a graph as in the statement of the lemma.
\end{proof}

For a multi-graph $\Delta$ we denote by $\doverline{\Delta}$ the ordinary graph obtained from $\Delta$ by replacing all multiple edges by single ones.

If two edges  of $\Delta$ are on the same vertices, then in its line graphs $L_1(\Delta)$
and $L_{\geq 1}(\Delta)$ they are in relation $\equiv$ and $\bowtie$, respectively.
So, if $v,w$ are two  vertices of $\Delta$ on an edge $e$, then 
identifying the  vertex $\{v,w\}$ of $L(\underline{\Delta})$
with the  $\equiv$- or $\bowtie$-equivalence class of $e$ in  $L_1(\Delta)$ or $L_{\geq 1}(\Delta)$,
respectively, yields the following:

\begin{lemma}\label{quotient of line graph}
Let $\Delta$ be a multi-graph.
Then $\overline{L_1(\Delta)}=\overline{L(\underline{\Delta})}$
and $\widehat{L_{\geq 1}(\Delta)}=\widehat{L(\underline{\Delta})}$.
\end{lemma}

A multi-graph $\Delta$ that does not contain any (not necessarily induced) subgraph $\Delta_0$
 
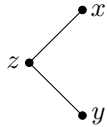
\begin{center}
 \begin{tikzpicture}[scale=0.7]
 \filldraw[color=black]
(0,0) circle [radius=2pt]
(1,1) circle [radius=2pt]
(1,-1) circle [radius=2pt];

\draw (1,1)--(0,0)--(1,-1);

\draw (1.3,1) node {$x$};
\draw (1.3,-1) node {$y$};
\draw (-0.3,0) node {$z$};

 \end{tikzpicture}
 \end{center}
 such that the vertices $x$ and $y$ have no other neighbors outside $\Delta_0$,
 is called \emph{$\Delta_0$-free}.
 
\begin{corollary}\label{cor}
Let $\Delta$ be  a connected multi-graph.
Then we have:

\begin{enumerate}
\item  
$\overline{L_1(\Delta)}={L(\underline{\Delta})}$,
or $\Delta$ has $4$ vertices.
\item
$\widehat{L_{\geq 1}(\Delta)}=L(\underline{\Delta})$,
or $\Delta$ is not $\Delta_0$-free.
\end{enumerate}
\end{corollary}

\begin{proof}
This follows from  \cref{except4} and \cref{quotient of line graph}   applied to
$\Gamma=L(\underline{\Delta})$.
\end{proof}

\begin{proposition}\label{iso}\label{4vertices}
Let $\Delta$ and $\Delta'$ be two connected multi-graphs.

\begin{enumerate}
\item Suppose $\phi$ 
is an isomorphism between $L_1(\Delta)$ and $L_1(\Delta')$.
Assume $\Delta$ and $\Delta'$ do not contain $4$ vertices.
Then  there is a unique isomorphism between $\Delta$ and $\Delta'$ inducing $\phi$.

\item 
Suppose $\phi$ 
is an isomorphism between $L_{\geq 1}(\Delta)$ and $L_{\geq 1}(\Delta')$.
Assume $\doverline \Delta$ and $\doverline \Delta'$ are $\Delta_0$-free.
Then  there is a unique isomorphism between $\Delta$ and $\Delta'$ inducing $\phi$.
\end{enumerate}
\end{proposition}

\begin{proof}
Suppose we are in case (i) or (ii).
By \cref{cor} and the assumptions,  $\overline{L_1(\Delta)}=L(\doverline{\Delta})$  and 
$\overline{L_1(\Delta')}=L(\doverline{\Delta'})$ are isomorphic in case (i), and 
$\widehat{L_1(\Delta)}=L(\doverline{\Delta})$  and 
$\widehat{L_1(\Delta')}=L(\doverline{\Delta'})$ are isomorphic in case (ii).

Then $\phi$ induces   a unique isomorphism
$\overline\phi$ between $\overline{L_1(\Delta)}$ and $\overline{L_1(\Delta')}$ in case (i)
and $\widehat\phi$ between $\widehat{L_1(\Delta)}$ and $\widehat{L_1(\Delta')}$ in case (ii).
Now, by Whitney's \cref{whitney} there is a unique isomorphism $\overline{\underline\phi}$ from  
$\doverline \Delta$ to $\doverline \Delta'$, which induces $\overline\phi$ between 
$\overline{L_1(\Delta)}=L(\doverline{\Delta})$ and $\overline{L_1(\Delta')}=L(\doverline{\Delta'})$.
As $K_{1,3}$ is not $\Delta_0$-free,
Whitney's \cref{whitney} also implies that there is a unique isomorphism $\widehat{\underline\phi}$ from  
$\doverline \Delta$ to $\doverline \Delta'$, which induces $\widehat\phi$ between 
$\widehat{L_{\geq 1}(\Delta)}=L(\doverline{\Delta})$ and $\widehat{L_{\geq 1}(\Delta')}=L(\doverline{\Delta'})$.

By the uniqueness of $\overline{\doverline\phi}$ and $\widehat{\doverline\phi}$
we find that the pairs 
$(\overline{\doverline\phi},\phi)$ and $(\widehat{\doverline\phi},\phi)$ are the unique isomorphisms
between $\Delta$ and $\Delta'$ we are looking for in case (i) and (ii), respectively.
\end{proof}

It remains to consider the cases excluded in the above \cref{iso}.
These are handled in the following result.

\begin{proposition}\label{nonisothm}
Let $\Delta$ be a connected  multi-graph.
\begin{enumerate}
\item
There is, up to isomorphism, a unique connected multi-graph $\Delta'$, not on $4$ vertices,
with $L_1(\Delta')$ isomorphic to $L_1(\Delta)$.
\item
There is, up to isomorphism, a unique
$\Delta_0$-free multi-graph $\Delta'$ with
$L_{\geq 1}(\Delta')$ and $L_{\geq 1}(\Delta)$ isomorphic.

\end{enumerate}

\end{proposition}

\begin{figure}
 \begin{center}
 \begin{tikzpicture}[scale=0.7]
 
 \draw[color=black,dashed]
(0,0) circle [radius=2];
\filldraw
(-2,0) circle [radius=2pt]
(-3,1) circle [radius=2pt]
(-3,0) circle [radius=2pt]
(-3,-1) circle [radius=2pt];
\filldraw
(2,0) circle [radius=2pt]
(3,0.5) circle [radius=2pt]
(3,-0.5) circle [radius=2pt];

\draw (-3,0)--(-2,0);
\draw (-3,1)--(-2,0)--(-3,-1);
\draw (3,0.5)--(2,0)--(3,-0.5)--(3,0.5);

\draw (-2,0)--(-1.5,0.5);
\draw (-2,0)--(-1.5,0);
\draw (-2,0)--(-1.5,-0.5);

\draw (2,0)--(1.5,0.5);
\draw (2,0)--(1.5,0);
\draw (2,0)--(1.5,-0.5);

 \end{tikzpicture}\hspace{.5cm}
 \begin{tikzpicture}[scale=0.7]
 
 \draw[color=black,dashed]
(0,0) circle [radius=2];
\filldraw
(-2,0) circle [radius=2pt]

(-3,0) circle [radius=2pt];
\filldraw
(2,0) circle [radius=2pt]
(3,0) circle [radius=2pt]
(4,0) circle [radius=2pt];

\draw (-3,0)--(-2,0);
\draw (-3,0.1)--(-2,0.1);
\draw (-3,-0.1)--(-2,-0.1);

\draw (2,0.1)--(3,0.1);
\draw (2,-0.1)--(3,-0.1);
\draw (3,0)--(4,0);

\draw (-2,0)--(-1.5,0.5);
\draw (-2,0)--(-1.5,0);
\draw (-2,0)--(-1.5,-0.5);

\draw (2,0)--(1.5,0.5);
\draw (2,0)--(1.5,0);
\draw (2,0)--(1.5,-0.5);

 \end{tikzpicture}
 \end{center}
 \caption{Two multi-graphs with the same $\geq 1$-line graph}
 \label{fig1}
 \end{figure}
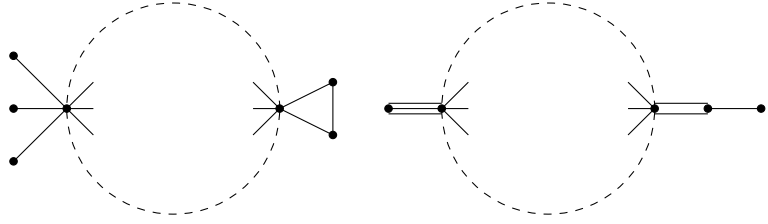

\begin{proof}
We first consider case (i).
By  \cref{iso}(i), we only have to consider the case that  $\Delta$
contains four vertices.

Then $\Delta$ is the multi-graph
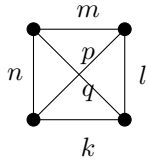
\begin{center}
\begin{tikzpicture}[scale=1.2]
\filldraw[color=black]
(0,0) circle [radius=2pt]
(1,0) circle [radius=2pt]
(1,1) circle [radius=2pt]
(0,1) circle [radius=2pt];

\draw (0,0)--(1,0)--(1,1)--(0,1)--(0,0);
\draw (0,0)--(1,1);
\draw (1,0)--(0,1);

\draw (0.6,-0.3) node {$k$};
\draw (0.6,1.2) node {$m$};
\draw (1.2,0.5) node {$l$};
\draw (-0.2,0.5) node {$n$};
\draw (0.6,0.7) node {$p$};
\draw (0.6,0.3) node {$q$};

\end{tikzpicture}
\end{center}
with $k,l,m$ cardinal numbers at least $1$ and $n,p,q$ cardinals at least $0$.

Let  $\Delta'$ be the following multi-graph:
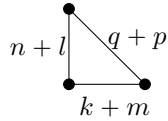
\begin{center}
\begin{tikzpicture}
\filldraw[color=black]
(0,0) circle [radius=2pt]
(1,0) circle [radius=2pt]
(0,1) circle [radius=2pt];

\draw (0,0)--(1,0)--(0,1)--(0,0);

\draw (0.6,-0.3) node {$k+m$};

\draw (-0.4,0.5) node {$n+l$};

\draw (0.9,0.6) node {$q+p$};

\end{tikzpicture}
\end{center}

Then $L_1(\Delta')$ is isomorphic to $L_1(\Delta)$
and by \cref{iso} it is the unique such graph which is not on $4$ vertices.

Now consider (ii).
We show that there does exist a $\Delta_0$-free root  multi-graph $\Delta'$ with
$L_{\geq 1}(\Delta')$ and $L_{\geq 1}(\Delta)$
isomorphic.

Suppose $\Delta$ does contain  subgraphs $\Delta_0$,
then we can consider the graph $\Delta'$ which we obtain from  
 $\Delta$  by identifying all the vertices in $[x]$ if $x$ and $y$ of $\Delta_0$ are not adjacent in $\Delta$ or ,if  $x$ and $y$  are not adjacent in $\Delta$,  by identifying $x$ and $y$ (so all vertices in $\langle x\rangle$) 
 and adding a new vertex $y'$ with for each edge 
 between $x$ and $y$ an edge between $\langle x\rangle$ and $y'$. See \cref{fig1}.
 
 By \cref{iso}(ii), this graph $\Delta'$ is then the unique $\Delta_0$-free root graph of $\Gamma$.
\end{proof}

 \cref{mainthm} and \cref{mainthm2}  follow from \cref{iso} and \cref{nonisothm}.

 \section{Find the root graph of a line graph}
\label{sect:alg}

Let $\Gamma$ be a finite graph.
In \cite{bermond}  it has been shown that $\Gamma$ is the $\geq 1$-line graph of a multi-graph $\Delta$
if and only if $\Gamma$ does not contain one of  $7$ graphs on at most $6$ vertices as induced subgraph.
A similar result for $1$-line graphs has been obtained in \cite{rootgraph}, where  a complete list of $33$ forbidden graphs, all on $6$ vertices, has been found.
This implies that there is a polynomial time algorithm to determine whether
$\Gamma$ is indeed a line graph.

The results of the previous section imply that, in  case $\Gamma$ is indeed a line graph, 
there is a unique root graph $\Delta$, not on four vertices if $\Gamma$ is a $1$-line graph,
or $\Delta_0$-free, if $\Gamma$ is a $\geq 1$-line graph.

The following algorithm has as on input of a finite graph $\Gamma$ as output this unique root graph $\Delta$ if $\Gamma$ is a $1$- or $\geq 1$-line graph, or returns the message
that the graph $\Gamma$ is not a line graph.

Correctness of the algorithm follows easily from the results of the previous section.

\begin{enumerate}
\item Find the false (or true) twin decomposition of $\Gamma$.
A decomposition into the $\equiv$-equivalence (or $\bowtie$-) classes can be found by a standard partition refinement.
See for example Algorithm 2 in \cite{partition}.
\item Determine $\overline \Gamma$ (or $\widehat{\Gamma}$). 
Take in each equivalence class a unique vertex and form the induced subgraph on these vertices.
\item Find a root graph $\underline \Delta$ for the graph $\overline \Gamma$ (or $\widehat{\Gamma}$), or return the message
that the graph $\Gamma$ is not a line graph. Several efficient algorithms exist. See for example \cite{alg3,alg1,alg2,alg4}.
\item Determine and output $\Delta$. Each vertex $v$ of $\Gamma$ 
determines a unique edge of $\overline \Gamma$ (or $\widehat{\Gamma}$) and hence of  $\underline \Delta$.
So,  we can assign to $v$ an edge of $\Delta$ on the same vertices of the corresponding edge of $\underline \Delta$.
\end{enumerate}

The above algorithm can of course also be used to determine a root graph of a generalized line graph.
Indeed, after step (i) one can check whether all $\equiv$-classes have size at most $2$,
in step (iv)  
whether for two vertices of $\Delta$ which are on two common edges, 
one is on no other edges.
In \cite{simic}  a different  algorithm is discussed to find a root graph of a generalized line graph.

Finally we notice that each of the steps in the above algorithm requires at most $O(|V|+|E|)$ computations. So, the complexity of the algorithm is $O(|V|+|E|)$.

\bibliographystyle{plain}

\bibliography{whitney.bib}

\vspace{2cm}

\parindent=0pt
Hans Cuypers\\
Department of Computer Science and Mathematics\\
Eindhoven University of Technology\\
P.O. Box 513 5600 MB, Eindhoven\\
The Netherlands\\
email: f.g.m.t.cuypers@tue.nl

 \end{document}